\documentclass{amsart}
\input{decls}
\UseAllTwocells
\title{Constructing symmetric monoidal bicategories}
\author{Michael A.\ Shulman}
\mdef\cMod{\mathcal{M}\mathit{od}}
\mdef\cCat{\mathcal{C}\mathit{at}}
\mdef\cTwocat{2\text{-}\mathcal{C}\mathit{at}}
\mdef\cBicat{\mathcal{B}\mathit{icat}}
\mdef\lMod{\mathbb{M}\mathsf{od}}
\mdef\lnCob{n\mathbb{C}\mathsf{ob}}
\mdef\lProf{\mathbb{P}\mathsf{rof}}
\mdef\cDbl{\mathcal{D}\mathit{bl}}
\mdef\fchk{\check{f}}
\mdef\conj{\Yleft}
\mdef\Conj{\mathcal{C}\mathit{onj}}
\begin{document}
\maketitle

\begin{abstract}
  We present a method of constructing symmetric monoidal bicategories
  from symmetric monoidal double categories that satisfy a lifting
  condition.  Such symmetric monoidal double categories frequently
  occur in nature, so the method is widely applicable, though not
  universally so.
\end{abstract}

\section{Introduction}
\label{sec:introduction}

Symmetric monoidal bicategories are important in many contexts.
However, the definition of even a monoidal bicategory
(see~\cite{gps:tricats,nick:tricats}), let alone a symmetric monoidal
one
(see~\cite{kv:2cat-zam,kv:bm2cat,bn:hda-i,ds:monbi-hopfagbd,crans:centers,mccrudden:bal-coalgb,gurski:brmonbicat}),
is quite imposing, and time-consuming to verify in any example.  In
this paper we describe a method for constructing symmetric monoidal
bicategories which is hardly more difficult than constructing a pair
of ordinary symmetric monoidal categories.  While not universally
applicable, this method applies in many cases of interest.  This idea
has often been implicitly used in particular cases, such as
bicategories of enriched profunctors, but to my knowledge the first
general statement was claimed in~\cite[Appendix B]{shulman:frbi}.  Our
purpose here is to work out the details, independently
of~\cite{shulman:frbi}.

\begin{rmk}
  Another approach to working out the details of this statement, from
  a different perspective, can be found
  in~\cite[\S5]{gg:ldstr-tricat}.  The two approaches contain
  basically the same content and results, although the authors
  of~\cite{gg:ldstr-tricat} work with ``locally-double bicategories''
  rather than monoidal double categories or 2x1-categories (see
  below).  They also don't treat the symmetric case, but as we will
  see, that is a fairly easy extension once the theory is in place.
  Thus, this note really presents nothing very new, only a
  self-contained and (hopefully) convenient treatment of the
  particular case of interest.
\end{rmk}

The method relies on the fact that in many bicategories, the 1-cells
are not the most fundamental notion of `morphism' between the objects.
For instance, in the bicategory \cMod\ of rings, bimodules, and
bimodule maps, the more fundamental notion of morphism between objects
is a ring homomorphism. The addition of these extra morphisms promotes
a bicategory to a \emph{double category}, or a category internal to
\cCat.  The extra morphisms are usually stricter than the 1-cells in
the bicategory and easier to deal with for coherence questions; in
many cases it is quite easy to show that we have a \emph{symmetric
  monoidal double category}.  The central observation is that in most
cases (when the double category is `fibrant') we can then `lift' this
symmetric monoidal structure to the original bicategory.  That is, we
prove the following theorem:

\begin{thm}\label{thm:mondbl-monbi-intro}
  If \lD\ is a fibrant monoidal double category, then its underlying
  bicategory $\cH(\lD)$ is a monoidal bicategory.  If \lD\ is braided
  or symmetric, so is $\cH(\lD)$.
\end{thm}

There is a good case to be made, however (see~\cite{shulman:frbi})
that often the extra morphisms should \emph{not} be discarded.  From
this point of view, in many cases symmetric monoidal bicategories are
a red herring, and really we should be studying symmetric monoidal
double categories.  This is also true in higher dimensions; for
instance, Chris Douglas~\cite{douglas:tfttalk} has suggested that
instead of tricategories we are usually interested in bicategories
internal to \cCat\ or categories internal to \cTwocat.  In most such
cases arising in practice, we can again `lift' the coherence to give a
tricategory after discarding the additional structure.

We propose the generic term \textbf{$(n\times k)$-category}
(pronounced ``$n$-by-$k$-category'') for an $n$-category internal to
$k$-categories, a structure which has $(n+1)(k+1)$ different types of
cells or morphisms arranged in an $(n+1)$ by $(k+1)$ grid.  Thus
double categories may be called \textbf{1x1-categories}, while in
place of tricategories we may consider 2x1-categories and
1x2-categories.  Any $(n\times k)$-category which satisfies a suitable
lifting property should have an underlying $(n+k)$-category, but
clearly as $n$ and $k$ grow an increasing amount of structure is
discarded in this process.


However, even for those of the opinion that $(n\times k)$-categories
are fundamental (such as the author), sometimes it really is the
underlying $(n+k)$-category that one cares about.  This is
particularly the case in the study of topological field theory, since
the Baez-Dolan cobordism hypothesis asserts a universal property of
the $(n+1)$-category of cobordisms which is not shared by the
$(n\times 1)$-category from which it is naturally constructed
(see~\cite{lurie:tft}).  Thus, regardless of one's philosophical bent,
results such as \autoref{thm:mondbl-monbi-intro} are of interest.

Proceeding to the contents of this paper, in
\S\ref{sec:symm-mono-double} we review the definition of symmetric
monoidal double categories, and in \S\ref{sec:comp-conj} we recall the
notions of `companion' and `conjoint' whose presence supplies the
necessary lifting property, which we call being \emph{fibrant}.  Then
in \S\ref{sec:1x1-to-bicat} we describe a functor from fibrant double
categories to bicategories, and in \S\ref{sec:constr-symm-mono} we
show that it preserves monoidal, braided, and symmetric structures.

I would like to thank Peter May, Tom Fiore, Stephan Stolz, Chris
Douglas, and Nick Gurski for helpful discussions and comments.

\section{Symmetric monoidal double categories}
\label{sec:symm-mono-double}

In this section, we introduce basic notions of double categories.
Double categories go back originally to Ehresmann
in~\cite{ehresmann:cat-str}; a brief introduction can be found
in~\cite{ks:r2cats}.  Other references
include~\cite{multi_funct_i,gp:double-limits,gp:double-adjoints}.

\begin{defn}
  A \textbf{(pseudo) double category} \lD\ consists of a `category of
  objects' $\lD_0$ and a `category of arrows' $\lD_1$, with structure
  functors
  \begin{align*}
    U&\maps \lD_0\to \lD_1\\
    S,T&\maps \lD_1\rightrightarrows \lD_0\\
    \odot&\maps \lD_1\times_{\lD_0}\lD_1\to \lD_1
  \end{align*}
  (where the pullback is over
  $\lD_1\too[T]\lD_0\overset{S}{\longleftarrow} \lD_1$) such that
  \begin{align*}
    S(U_A) &= A\\
    T(U_A) &= A\\
    S(M\odot N) &= SN\\
    T(M\odot N) &= TM
  \end{align*}
  equipped with natural isomorphisms
  \begin{align*}
    \fa &: (M\odot N) \odot P \too[\iso] M \odot (N \odot P)\\
    \fl &: U_B \odot M \too[\iso] M\\
    \fr &: M \odot U_A \too[\iso] M
  \end{align*}
  such that $S(\fa)$, $T(\fa)$, $S(\fl)$, $T(\fl)$, $S(\fr)$, and
  $T(\fr)$ are all identities, and such that the standard coherence
  axioms for a monoidal category or bicategory (such as Mac Lane's
  pentagon; see~\cite{maclane}) are satisfied.
\end{defn}

Just as a bicategory can be thought of as a category weakly
\emph{enriched} over \cCat, a pseudo double category can be thought of
as a category weakly \emph{internal} to \cCat.  Since these are the
kind of double category of most interest to us, we will usually drop
the adjective ``pseudo.''

We call the objects of $\lD_0$ \textbf{objects} or \textbf{0-cells},
and we call the morphisms of $\lD_0$ \textbf{(vertical) 1-morphisms}
and write them as $f\maps A\to B$.  We call the objects of $\lD_1$
\textbf{(horizontal) 1-cells}; if $M$ is a 1-cell with $S(M)=A$ and
$T(M)=B$, we write $M\maps A\hto B$.  We call a morphism $\alpha\maps
M\to N$ of $\lD_1$ with $S(\alpha)=f$ and $T(\alpha)=g$ a
\textbf{2-morphism} and draw it as follows:
\begin{equation}\label{eq:square}
  \xymatrix@-.5pc{
    A \ar[r]|{|}^{M}  \ar[d]_f \ar@{}[dr]|{\Downarrow\alpha}&
    B\ar[d]^g\\
    C \ar[r]|{|}_N & D
  }.
\end{equation}
Note that we distinguish between \emph{1-morphisms}, which we draw
vertically, and \emph{1-cells}, which we draw horizontally.  In
traditional double-category terminology these are both referred to
with the same word (be it ``cell'' or ``morphism'' or ``arrow''), the
distinction being made only by the adjectives ``vertical'' and
``horizontal.''  Our terminology is more concise, and allows for
flexibility in the drawing of pictures without a corresponding change
in names (some authors prefer to draw their double categories
transposed from ours).

We write the composition of vertical 1-morphisms $A\too[f] B\too[g] C$
and the vertical composition of 2-morphisms $M\too[\alpha]
N\too[\beta] P$ as $g\circ f$ and $\beta\circ\alpha$, or sometimes
just $gf$ and $\beta\alpha$.  We write the horizontal composition of
1-cells $A\xhto{M} B \xhto{N} C$ as $A\xhto{N\odot M} C$ and that of
2-morphisms
\[\vcenter{\xymatrix{ \ar[r]|-@{|}^-{} \ar[d] \ar@{}[dr]|{\Downarrow\alpha} &
     \ar[r]|-@{|}^-{} \ar[d] \ar@{}[dr]|{\Downarrow\beta} &\ar[d]\\
  \ar[r]|-@{|}_-{} & \ar[r]|-@{|}_-{} & }}\]
as
\[\vcenter{\xymatrix@C=4pc{ \ar[r]|-@{|}^-{} \ar[d] \ar@{}[dr]|{\Downarrow\;\be\odot\al} &  \ar[d]\\
  \ar[r]|-@{|}_-{} & }}\]
The two different compositions of 2-morphisms obey an interchange law,
by the functoriality of $\odot$:
\[(M_1\odot M_2) \circ (N_1\odot N_2) = (M_1\circ N_1)\odot (M_2\circ N_2).
\]
Every object $A$ has a vertical identity $1_A$ and a horizontal unit
$U_A$, every 1-cell $M$ has an identity 2-morphism $1_M$, every
vertical 1-morphism $f$ has a horizontal unit 2-morphism $U_f$, and we
have $1_{U_A} = U_{1_A}$ (by the functoriality of $U$).

Note that the vertical composition $\circ$ is strictly associative and
unital, while the horizontal one $\odot$ is only weakly so.  This is
the case in most of the examples we have in mind.  It is possible to
define double categories that are weak in both directions (see, for
instance,~\cite{verity:base-change}), but this introduces much more
complication and is usually unnecessary.

\begin{rmk}\label{rmk:monglob}
  In general, an $(n\times 1)$-category consists of 1-categories
  $\lD_i$ for $0\le i\le n$, together with source, target, unit, and
  composition functors and coherence isomorphisms.  We refer to the
  objects of $\lD_i$ as \textbf{$i$-cells} and to the morphisms of
  $\lD_i$ as \textbf{morphisms of $i$-cells} or \textbf{(vertical)
    $(i+1)$-morphisms}.  A formal definition can be found
  in~\cite{batanin:monglob} under the name \emph{monoidal $n$-globular
    category}.
\end{rmk}


A 2-morphism~\eqref{eq:square} where $f$ and $g$ are identities (such
as the constraint isomorphisms $\fa,\fl,\fr$) is called
\textbf{globular}.  Every double category \lD\ has a
\textbf{horizontal bicategory} $\cH(\lD)$ consisting of the objects,
1-cells, and globular 2-morphisms.  Conversely, many naturally
occurring bicategories are actually the horizontal bicategory of some
naturally ocurring double category.  Here are just a few examples.

\begin{eg}
  The double category \lMod\ has as objects rings, as 1-morphisms ring
  homomorphisms, as 1-cells bimodules, and as 2-morphisms equivariant
  bimodule maps.  Its horizontal bicategory $\cMod = \cH(\lMod)$ is
  the usual bicategory of rings and bimodules.
\end{eg}

\begin{eg}
  The double category \lnCob\ has as objects closed $n$-manifolds, as
  1-morphisms diffeomorphisms, as 1-cells cobordisms, and as
  2-morphisms diffeomorphisms between cobordisms.  Again $\cH(\lnCob)$
  is the usual bicategory of cobordisms.
\end{eg}

\begin{eg}
  The double category \lProf\ has as objects categories, as
  1-morphisms functors, as 1-cells \emph{profunctors} (a profunctor
  $A\hto B$ is a functor $B\op\times A\to \mathbf{Set}$), and as
  2-morphisms natural transformations.  Bicategories such as
  $\cH(\lProf)$ are commonly encountered in category theory,
  especially the enriched versions.
\end{eg}


As opposed to bicategories, which naturally form a tricategory, double
categories naturally form a \emph{2-category}, a much simpler object.

\begin{defn}
  Let \lD\ and \lE\ be double categories.  A \textbf{(pseudo double)
    functor} $F\maps \lD\to \lE$ consists of the following.
  \begin{itemize}
  \item Functors $F_0\maps \lD_0 \to \lE_0$ and $F_1\maps \lD_1 \to
    \lE_1$ such that $S\circ F_1 = F_0\circ S$ and $T\circ F_1 =
    F_0\circ T$.
  \item Natural transformations $F_\odot\maps F_1M \odot F_1N \to
    F_1(M\odot N)$ and $F_U\maps U_{F_0 A} \to F_1(U_A)$, whose
    components are globular isomorphisms, and which satisfy the usual
    coherence axioms for a monoidal functor or pseudofunctor
    (see~\cite[\S{}XI.2]{maclane}).
  \end{itemize}
\end{defn}

\begin{defn}\label{thm:dbl-transf}
  A \textbf{(vertical) transformation} between two functors $\alpha:
  F\to G:\lD\to\lE$ consists of natural transformations $\alpha_0\maps
  F_0\to G_0$ and $\alpha_1\maps F_1\to G_1$ (both usually written as
  $\alpha$), such that $S(\alpha_{M}) = \alpha_{SM}$ and
  $T(\alpha_{M}) = \alpha_{TM}$, and such that
  \[\vcenter{\xymatrix@-.5pc{
      FA \ar@{=}[d] \ar[r]|{|}^{FM}
      \ar@{}[drr]|{\Downarrow F_\odot} &
      FB \ar[r]|{|}^{FN} &
      FC \ar@{=}[d]\\
      FA \ar[rr]|{F(N\odot M)} \ar[d]_{\alpha_A}
      \ar@{}[drr]|{\Downarrow \alpha_{N\odot M}} &&
      FC \ar[d]^{\alpha_C}\\
      GA \ar[rr]|{|}_{G(N\odot M)} && GC
    }} =
  \vcenter{\xymatrix@-.5pc{
      FA \ar[d]_{\alpha_A} \ar@{}[dr]|{\Downarrow \alpha_M} \ar[r]|{|}^{FM} &
      FB \ar[d]|{\alpha_B} \ar@{}[dr]|{\Downarrow \alpha_N} \ar[r]|{|}^{FN} &
      FC \ar[d]^{\alpha_C}\\
      GA \ar@{=}[d] \ar[r]|{|}_{GM} \ar@{}[drr]|{\Downarrow G_\odot} &
      GB \ar[r]|{|}_{GN} &
      GC \ar@{=}[d]\\
      GA \ar[rr]|{|}_{G(N\odot M)} && GC
    }}\]
  for all 1-cells $M\colon A\hto B$ and $N\colon B\hto C$, and
  \[\vcenter{\xymatrix@-.5pc{
      FA \ar[rr]|{|}^{U_{FA}} \ar@{=}[d]
      \ar@{}[drr]|{\Downarrow F_0} &&
      FA \ar@{=}[d]\\
      FA \ar[rr]|{F(U_A)} \ar[d]_{\alpha_A}
      \ar@{}[drr]|{\Downarrow \alpha_{U_A}} &&
      FA \ar[d]^{\alpha_A}\\
      GA \ar[rr]|{|}_{G(U_A)} && GA
    }} =
  \vcenter{\xymatrix@-.5pc{
      FA \ar[rr]|{|}^{U_{FA}} \ar[d]_{\alpha_A}
      \ar@{}[drr]|{\Downarrow U_{\alpha_A}} &&
      FA \ar[d]^{\alpha_A}\\
      GA \ar[rr]|{U_{GA}} \ar@{=}[d]
      \ar@{}[drr]|{\Downarrow F_0} &&
      GA \ar@{=}[d]\\
      GA \ar[rr]|{|}_{G(U_A)} && GA.
    }}\]
  for all objects $A$.
\end{defn}

We write \cDbl\ for the 2-category of double categories, functors, and
transformations, and $\mathbf{Dbl}$ for its underlying 1-category.
Note that a 2-cell $\al$ in \cDbl\ is an isomorphism just when each
$\al_A$, \emph{and} each $\al_M$, is invertible.

The 2-category \cDbl\ gives us an easy way to define what we mean by a
\emph{symmetric monoidal double category}.  In any 2-category with
finite products there is a notion of a \emph{pseudomonoid}, which
generalizes the notion of monoidal category in \cCat.  Specializing
this to \cDbl, we obtain the following.

\begin{defn}
  A \textbf{monoidal double category} is a double category equipped
  with functors $\ten\maps \lD\times\lD\to\lD$ and $I\maps * \to\lD$,
  and invertible transformations
  \begin{align*}
    \mathord{\otimes} \circ (\Id\times \mathord{\otimes})
    &\iso \mathord{\otimes} \circ (\mathord{\otimes} \times \Id)\\
    \mathord{\otimes} \circ (\Id\times I) &\iso \Id\\
    \mathord{\otimes} \circ (I\times \Id) &\iso \Id
  \end{align*}
  satisfying the usual axioms.  If it additionally has a braiding
  isomorphism
  \begin{align*}
    \mathord{\otimes} &\iso \mathord{\otimes} \circ \tau
  \end{align*}
  (where $\tau\maps \lD\times\lD\iso \lD\times\lD$ is the twist)
  satisfying the usual axioms, then it is \textbf{braided} or
  \textbf{symmetric}, according to whether or not the braiding is
  self-inverse.
\end{defn}

Unpacking this definition more explicitly, we see that a monoidal
double category is a double category together with the following
structure.
\begin{enumerate}
\item $\lD_0$ and $\lD_1$ are both monoidal categories.
\item If $I$ is the monoidal unit of $\lD_0$, then $U_I$ is the
  monoidal unit of $\lD_1$.\footnote{Actually, all the above
    definition requires is that $U_I$ is coherently \emph{isomorphic
      to} the monoidal unit of $\lD_1$, but we can always choose them
    to be equal without changing the rest of the structure.}
\item The functors $S$ and $T$ are strict monoidal, i.e.\ $S(M\ten N)
  = SM\ten SN$ and $T(M\ten N)=TM\ten TN$ and $S$ and $T$ also
  preserve the associativity and unit constraints.
\item We have globular isomorphisms
  \[\fx\maps (M_1\ten N_1)\odot (M_2\ten N_2)\too[\iso] (M_1\odot M_2)\ten (N_1\odot N_2)\]
  and
  \[\fu\maps U_{A\ten B} \too[\iso] (U_A \ten U_B)\]
  such that the following diagrams commute:
  \[\xymatrix{
    ((M_1\ten N_1)\odot (M_2\ten N_2)) \odot (M_3\ten N_3) \ar[r]\ar[d]
    & ((M_1\odot M_2)\ten (N_1\odot N_2)) \odot (M_3\ten N_3) \ar[d]\\
    (M_1\ten N_1)\odot ((M_2\ten N_2) \odot (M_3\ten N_3)) \ar[d] &
    ((M_1\odot M_2)\odot M_3) \ten ((N_1\odot N_2)\odot N_3) \ar[d]\\
    (M_1\ten N_1) \odot ((M_2\odot M_3) \ten (N_2\odot N_3))\ar[r] &
    (M_1\odot (M_2\odot M_3)) \ten (N_1\odot (N_2\odot N_3))}\]
  \[\xymatrix{(M\ten N) \odot U_{C\ten D} \ar[r]\ar[d] &
    (M\ten N)\odot (U_C\ten U_D) \ar[d]\\
    M\ten N\ar@{<-}[r] & (M\odot U_C) \ten (N\odot U_D)}\]
  \[\xymatrix{U_{A\ten B}\odot (M\ten N)  \ar[r]\ar[d] &
    (U_A\ten U_B)\odot (M\ten N) \ar[d]\\
    M\ten N\ar@{<-}[r] & (U_A \odot M) \ten (U_B\odot N)}\]
  (these arise from the constraint data for the pseudo double functor
  $\ten$).
\item The following diagrams commute, expressing that the
  associativity isomorphism for $\ten$ is a transformation of double
  categories.
  \[\xymatrix{
    ((M_1\ten N_1)\ten P_1) \odot ((M_2\ten N_2)\ten P_2) \ar[r]\ar[d] &
    (M_1\ten (N_1\ten P_1)) \odot (M_2\ten (N_2\ten P_2)) \ar[d]\\
    ((M_1\ten N_1) \odot (M_2\ten N_2)) \ten (P_1\odot P_2) \ar[d] &
    (M_1\odot M_2) \ten ((N_1\ten P_1)\odot (N_2\ten P_2))\ar[d] \\
    ((M_1\odot M_2) \ten(N_1\odot N_2)) \ten (P_1\odot P_2) \ar[r] &
    (M_1\odot M_2) \ten ((N_1\odot N_2)\ten (P_1\odot P_2))}\]
  \[\xymatrix{
    U_{(A\ten B)\ten C} \ar[r] \ar[d] & U_{A\ten (B\ten C)} \ar[d]\\
    U_{A\ten B} \ten U_C \ar[d] & U_A\ten U_{B\ten C}\ar[d]\\
    (U_A\ten U_B)\ten U_C \ar[r] & U_A\ten (U_B\ten U_C) }\]
\item The following diagrams commute, expressing that the unit
  isomorphisms for $\ten$ are transformations of double categories.
  \[\vcenter{\xymatrix{
      (M\ten U_I)\odot (N\ten U_I)\ar[r]\ar[d] &
      (M\odot N)\ten (U_I \odot U_I) \ar[d]\\
      M\odot N \ar@{<-}[r] &
      (M\odot N)\ten U_I }}\]
  \[\vcenter{\xymatrix{U_{A\ten I} \ar[r]\ar[dr] & U_A\ten U_I \ar[d]\\
       & U_A}}\]
  \[\vcenter{\xymatrix{
      (U_I\ten M)\odot (U_I\ten N)\ar[r]\ar[d] &
      (U_I \odot U_I) \ten (M\odot N) \ar[d]\\
      M\odot N \ar@{<-}[r] &
      U_I\ten (M\odot N) }}\]
  \[\vcenter{\xymatrix{U_{I\ten A} \ar[r]\ar[dr] & U_I\ten U_A \ar[d]\\
      & U_A}}\]
  \newcounter{mondbl}
  \setcounter{mondbl}{\value{enumi}}
\end{enumerate}
Similarly, a braided monoidal double category is a monoidal double
category with the following additional structure.
\begin{enumerate}\setcounter{enumi}{\value{mondbl}}
\item $\lD_0$ and $\lD_1$ are braided monoidal categories.
\item The functors $S$ and $T$ are strict braided monoidal (i.e.\ they
  preserve the braidings).
\item The following diagrams commute, expressing that the braiding is
  a transformation of double categories.
  \[\xymatrix{(M_1\odot M_2)\ten (N_1\odot N_2) \ar[r]^\fs\ar[d]_\fx &
    (N_1\odot N_2)\ten (M_1 \odot M_2)\ar[d]^\fx\\
    (M_1\ten N_1)\odot (M_2\ten N_2) \ar[r]_{\fs\odot \fs} &
    (N_1\ten M_1) \odot (N_2 \ten M_2)}
  \]
  \[\xymatrix{U_A \ten U_B \ar[r]^(0.55)\fu \ar[d]_\fs &
    U_{A\ten B} \ar[d]^{U_\fs}\\
    U_B\ten U_A \ar[r]_(0.55)\fu &
    U_{B\ten A}}.
  \]
  \setcounter{mondbl}{\value{enumi}}
\end{enumerate}
Finally, a symmetric monoidal double category is a braided one such that
\begin{enumerate}\setcounter{enumi}{\value{mondbl}}
\item $\lD_0$ and $\lD_1$ are in fact symmetric monoidal.
\end{enumerate}
While there are a fair number of coherence diagrams to verify, most of
them are fairly small, and in any given case most or all of them are
fairly obvious.  Thus, verifying that a given double category is
(braided or symmetric) monoidal is not a great deal of work.

\begin{eg}
  The examples \lMod, \lnCob, and \lProf\ are all easily seen to be
  symmetric monoidal under the tensor product of rings, disjoint union
  of manifolds, and cartesian product of categories, respectively.
\end{eg}

\begin{rmk}
  In a 2-category with finite products there is additionally the
  notion of a \emph{cartesian object}: one such that the diagonal
  $D\to D\times D$ and projection $D\to 1$ have right adjoints.  Any
  cartesian object is a symmetric pseudomonoid in a canonical way,
  just as any category with finite products is a monoidal category
  with its cartesian product.  Many of the ``cartesian bicategories''
  considered in~\cite{cw:cart-bicats-i,ckww:cartbicats-ii} are in
  fact the horizontal bicategory of some cartesian object in \cDbl,
  and inherit their monoidal structure in this way.
\end{rmk}

Two further general methods for constructing symmetric monoidal double
categories can be found in~\cite{shulman:frbi}.

\begin{rmk}
  The general yoga of internalization says that an $X$ internal to
  $Y$s internal to $Z$s is equivalent to a $Y$ internal to $X$s
  internal to $Z$s, but this is only strictly true when the
  internalizations are all strict.  We have defined a symmetric
  monoidal double category to be a (pseudo) symmetric monoid internal
  to (pseudo) categories internal to categories, but one could also
  consider a (pseudo) category internal to (pseudo) symmetric monoids
  internal to categories, i.e.\ a pseudo internal category in the
  2-category
  $\mathcal{S}\mathit{ym}\mathcal{M}\mathit{on}\mathcal{C}\mathit{at}$
  of symmetric monoidal categories and strong symmetric monoidal
  functors.  This would give \emph{almost} the same definition, except
  that $S$ and $T$ would only be strong monoidal (preserving $\ten$ up
  to isomorphism) rather than strict monoidal.  We prefer our
  definition, since $S$ and $T$ are strict monoidal in almost all
  examples, and keeping track of their constraints would be tedious.
\end{rmk}

Just as every bicategory is equivalent to a strict 2-category, it is
proven in~\cite{gp:double-limits} that every pseudo double category is
equivalent to a strict double category (one in which the associativity
and unit constraints for $\odot$ are identities).  Thus, from now on
we will usually omit to write these constraint isomorphisms (or
equivalently, implicitly strictify our double categories).  We
\emph{will} continue to write the constraint isomorphisms for the
monoidal structure $\ten$, since these are where the whole question
lies.

\section{Companions and conjoints}
\label{sec:comp-conj}

Suppose that \lD\ is a symmetric monoidal double category; when does
$\cH(\lD)$ become a symmetric monoidal bicategory?  It clearly has a
unit object $I$, and the pseudo double functor $\ten\maps
\lD\times\lD\to\lD$ clearly induces a functor $\ten\maps
\cH(\lD)\times\cH(\lD)\to\cH(\lD)$.  However, the problem is that the
constraint isomorphisms such as $A\ten (B\ten C)\iso (A\ten B)\ten C$
are \emph{vertical} 1-morphisms, which get discarded when we pass to
$\cH(\lD)$.  Thus, in order for $\cH(\lD)$ to inherit a symmetric
monoidal structure, we must have a way to make vertical 1-morphisms
into horizontal 1-cells.  Thus is the purpose of the following
definition.

\begin{defn}\label{def:companion}
  Let \lD\ be a double category and $f\maps A\to B$ a vertical
  1-morphism.  A \textbf{companion} of $f$ is a horizontal 1-cell
  $\fhat\maps A\hto B$ together with 2-morphisms
  \begin{equation*}
    \begin{array}{c}
      \xymatrix@-.5pc{
        \ar[r]|-@{|}^-{\fhat} \ar[d]_f \ar@{}[dr]|\Downarrow
        & \ar@{=}[d]\\
        \ar[r]|-@{|}_-{U_B} & }
    \end{array}\quad\text{and}\quad
    \begin{array}{c}
      \xymatrix@-.5pc{
        \ar[r]|-@{|}^-{U_A} \ar@{=}[d] \ar@{}[dr]|\Downarrow
        & \ar[d]^f\\
        \ar[r]|-@{|}_-{\fhat} & }
    \end{array}
  \end{equation*}
  such that the following equations hold.
  \begin{align}\label{eq:compeqn}
    \begin{array}{c}
      \xymatrix@-.5pc{
        \ar[r]|-@{|}^-{U_A} \ar@{=}[d] \ar@{}[dr]|\Downarrow
        & \ar[d]^f\\
        \ar[r]|-{\fhat} \ar[d]_f \ar@{}[dr]|\Downarrow
        & \ar@{=}[d]\\
        \ar[r]|-@{|}_-{U_B} & }
    \end{array} &= 
    \begin{array}{c}
      \xymatrix@-.5pc{ \ar[r]|-@{|}^-{U_A} \ar[d]_f
        \ar@{}[dr]|{\Downarrow U_f} &  \ar[d]^f\\
        \ar[r]|-@{|}_-{U_B} & }
    \end{array}
    &
    \begin{array}{c}
      \xymatrix@-.5pc{
        \ar[r]|-@{|}^-{U_A} \ar@{=}[d] \ar@{}[dr]|\Downarrow &
        \ar[r]|-@{|}^-{\fhat} \ar[d]_f \ar@{}[dr]|\Downarrow
        & \ar@{=}[d]\\
        \ar[r]|-@{|}_-{\fhat} &
        \ar[r]|-@{|}_-{U_B} &}
    \end{array} &=
    \begin{array}{c}
      \xymatrix@-.5pc{
        \ar[r]|-@{|}^-{\fhat} \ar@{=}[d] \ar@{}[dr]|{\Downarrow 1_{\fhat}}
        & \ar@{=}[d]\\
        \ar[r]|-@{|}_-{\fhat} & }
    \end{array}
  \end{align}
  A \textbf{conjoint} of $f$, denoted $\fchk\maps B\hto A$, is a
  companion of $f$ in the double category $\lD^{h\cdot\mathrm{op}}$
  obtained by reversing the horizontal 1-cells, but not the vertical
  1-morphisms, of \lD.
\end{defn}

\begin{rmk}
  We momentarily suspend our convention of pretending that our double
  categories are strict to mention that the second
  equation in~\eqref{eq:compeqn} actually requires an insertion of unit
  isomorphisms to make sense.
\end{rmk}

The form of this definition is due
to~\cite{gp:double-adjoints,dpp:spans}, but the ideas date back
to~\cite{bs:dblgpd-xedmod}; see
also~\cite{bm:dbl-thin-conn,fiore:pscat}.  In the terminology of these
references, a \emph{connection} on a double category is equivalent to
a strictly functorial choice of a companion for each vertical arrow.

\begin{defn}
  We say that a double category is \textbf{fibrant} if every vertical
  1-morphism has both a companion and a conjoint.
\end{defn}

\begin{rmk}
  In~\cite{shulman:frbi} fibrant double categories were called
  \emph{framed bicategories}.  However, the present terminology seems
  to generalize better to $(n\times k)$-categories, as well as
  avoiding a conflict with the \emph{framed bordisms} in topological
  field theory.
\end{rmk}


\begin{egs}
  \lMod, \lnCob, and \lProf\ are all fibrant.  In \lMod, the companion
  of a ring homomorphism $f\maps A\to B$ is $B$ regarded as an
  $A$-$B$-bimodule via $f$ on the left, and dually for its conjoint.
  In \lnCob, companions and conjoints are obtained by regarding a
  diffeomorphism as a cobordism.  And in \lProf, companions and
  conjoints are obtained by regarding a functor $f\maps A\to B$ as a
  `representable' profunctor $B(f-,-)$ or $B(-,f-)$.
\end{egs}

\begin{rmk}
  For an $(n\times 1)$-category (recall \autoref{rmk:monglob}), the
  lifting condition we should require is simply that each double
  category $\lD_{i+1} \toto \lD_i$, for $0\le i < n$, is fibrant.
\end{rmk}

The existence of companions and conjoints gives us a way to `lift'
vertical 1-morphisms to horizontal 1-cells.  What is even more crucial
for our applications, however, is that these liftings are unique up to
isomorphism, and that these isomorphisms are canonical and coherent.
This is the content of the following lemmas.  We state most of them
only for companions, but all have dual versions for conjoints.

\begin{lem}\label{thm:theta}
  Let $\fhat\maps A\hto B$ and $\fhat'\maps A\hto B$ be companions of
  $f$ (that is, each comes \emph{equipped with} 2-morphisms as in
  \autoref{def:companion}).  Then there is a unique globular isomorphism
  $\theta_{\fhat,\fhat'}\maps \fhat\too[\iso]\fhat'$ such that
  \begin{equation}\label{eq:comp-iso}
    \vcenter{\xymatrix@R=1.5pc{
        \ar[r]|-@{|}^-{U_A} \ar@{=}[d] \ar@{}[dr]|\Downarrow &  \ar[d]^f\\
        \ar[r]|-{\fhat} \ar@{=}[d] \ar@{}[dr]|{\Downarrow \theta_{\fhat,\fhat'}} &  \ar@{=}[d]\\
        \ar[r]|-{\fhat'} \ar[d]_f \ar@{}[dr]|\Downarrow &  \ar@{=}[d]\\
        \ar[r]|-@{|}_-{U_B} & }} \quad = \quad
    \vcenter{\xymatrix@-.5pc{ \ar[r]|-@{|}^-{U_A} \ar[d]_f
        \ar@{}[dr]|{\Downarrow U_f} &  \ar[d]^f\\
        \ar[r]|-@{|}_-{U_B} & .}}
  \end{equation}
\end{lem}
\begin{proof}
  Composing~\eqref{eq:comp-iso} on the left with
  $\vcenter{\xymatrix@-.5pc{ \ar[r]|-@{|}^-{U_A} \ar@{=}[d]
      \ar@{}[dr]|\Downarrow & \ar[d]^f\\ \ar[r]|-@{|}_-{\fhat'} & }}$
  and on the right with $\vcenter{\xymatrix@-.5pc{
      \ar[r]|-@{|}^-{\fhat} \ar[d]_f \ar@{}[dr]|\Downarrow &
      \ar@{=}[d]\\ \ar[r]|-@{|}_-{U_B} & }}$, and using the second
  equation~\eqref{eq:compeqn}, we see that if~\eqref{eq:comp-iso} is
  satisfied then $\theta_{\fhat,\fhat'}$ must be the composite
  \begin{equation}
    \vcenter{\xymatrix@-.5pc{
        \ar[r]|-@{|}^-{U_A} \ar@{=}[d] \ar@{}[dr]|\Downarrow &
        \ar[r]|-@{|}^-{\fhat} \ar[d]|f \ar@{}[dr]|\Downarrow
        & \ar@{=}[d]\\
        \ar[r]|-@{|}_-{\fhat'} &
        \ar[r]|-@{|}_-{U_B} &}}\label{eq:theta}
  \end{equation}
  Two applications of the first equation~\eqref{eq:compeqn} shows that
  this indeed satisfies~\eqref{eq:comp-iso}.  As for its being an
  isomorphism, we have the dual composite $\theta_{\fhat',\fhat'}$:
  \[\vcenter{\xymatrix@-.5pc{
      \ar[r]|-@{|}^-{U_A} \ar@{=}[d] \ar@{}[dr]|\Downarrow &
      \ar[r]|-@{|}^-{\fhat'} \ar[d]_f \ar@{}[dr]|\Downarrow
      & \ar@{=}[d]\\
      \ar[r]|-@{|}_-{\fhat} &
      \ar[r]|-@{|}_-{U_B} &}}\]
  which we verify is an inverse using~\eqref{eq:compeqn}:
  \[\vcenter{\xymatrix@-.5pc{
      \ar[r]|-@{|}^{U_A}\ar@{=}[d] \ar@{}[dr]|{=} &
      \ar[r]|-@{|}^{U_A}\ar@{=}[d] \ar@{}[dr]|{\Downarrow} &
      \ar[r]|-@{|}^{\fhat}\ar[d]|f \ar@{}[dr]|{\Downarrow} &
      \ar@{=}[d]\\
      \ar[r]|{U_A}\ar@{=}[d] \ar@{}[dr]|{\Downarrow} &
      \ar[r]|{\fhat'}\ar[d]|f \ar@{}[dr]|{\Downarrow} &
      \ar[r]|{U_B}\ar@{=}[d] \ar@{}[dr]|{=} &
      \ar@{=}[d]\\
      \ar[r]|-@{|}_{\fhat} &
      \ar[r]|-@{|}_{U_B} &
      \ar[r]|-@{|}_{U_B} &
    }} \;=\;
  \vcenter{\xymatrix@-.5pc{
      \ar[r]|-@{|}^-{U_A} \ar@{=}[d] \ar@{}[dr]|\Downarrow &
      \ar[r]|-@{|}^-{\fhat} \ar[d]_f \ar@{}[dr]|\Downarrow
      & \ar@{=}[d]\\
      \ar[r]|-@{|}_-{\fhat} &
      \ar[r]|-@{|}_-{U_B} &}} \;=\;
  \vcenter{\xymatrix@-.5pc{
      \ar[r]|-@{|}^-{\fhat} \ar@{=}[d] \ar@{}[dr]|{\Downarrow 1_{\fhat}}
      & \ar@{=}[d]\\
      \ar[r]|-@{|}_-{\fhat} & }}\]
  (and dually).
\end{proof}

\begin{lem}\label{thm:theta-id}
  For any companion \fhat\ of $f$ we have $\theta_{\fhat,\fhat}=1_{\fhat}$.
\end{lem}
\begin{proof}
  This is the second equation~\eqref{eq:compeqn}.
\end{proof}

\begin{lem}\label{thm:theta-compose-vert}
  Suppose that $f$ has three companions $\fhat$, $\fhat'$, and
  $\fhat''$.  Then $\theta_{\fhat,\fhat''} = \theta_{\fhat',\fhat''}
  \circ\theta_{\fhat,\fhat'}$.
\end{lem}
\begin{proof}
  By definition, we have
  \[\theta_{\fhat',\fhat''} \circ\theta_{\fhat,\fhat'} =\;
  \vcenter{\xymatrix@-.5pc{
      \ar[r]|-@{|}^{U_A}\ar@{=}[d] \ar@{}[dr]|{=} &
      \ar[r]|-@{|}^{U_A}\ar@{=}[d] \ar@{}[dr]|{\Downarrow} &
      \ar[r]|-@{|}^{\fhat}\ar[d]|f \ar@{}[dr]|{\Downarrow} &
      \ar@{=}[d]\\
      \ar[r]|{U_A}\ar@{=}[d] \ar@{}[dr]|{\Downarrow} &
      \ar[r]|{\fhat'}\ar[d]|f \ar@{}[dr]|{\Downarrow} &
      \ar[r]|{U_B}\ar@{=}[d] \ar@{}[dr]|{=} &
      \ar@{=}[d]\\
      \ar[r]|-@{|}_{\fhat''} &
      \ar[r]|-@{|}_{U_B} &
      \ar[r]|-@{|}_{U_B} &
    }} \;=\;
  \vcenter{\xymatrix@-.5pc{
      \ar[r]|-@{|}^-{U_A} \ar@{=}[d] \ar@{}[dr]|\Downarrow &
      \ar[r]|-@{|}^-{\fhat} \ar[d]_f \ar@{}[dr]|\Downarrow
      & \ar@{=}[d]\\
      \ar[r]|-@{|}_-{\fhat''} &
      \ar[r]|-@{|}_-{U_B} &}} \;=
  \theta_{\fhat,\fhat''}\]
  as desired.
\end{proof}

\begin{lem}\label{thm:comp-unit}
  $U_A\maps A\hto A$ is always a companion of $1_A\maps A\to A$ in a
  canonical way.
\end{lem}
\begin{proof}
  We take both defining 2-morphisms to be
  $1_{U_A}$; the truth of~\eqref{eq:compeqn} is evident.
\end{proof}

\begin{lem}\label{thm:comp-compose}
  Suppose that $f\maps A\to B$ has a companion \fhat\ and $g\maps B\to
  C$ has a companion \ghat.  Then $\ghat\odot\fhat$ is a companion of
  $gf$.
\end{lem}
\begin{proof}
  We take the defining 2-morphisms to be the composites
  \[\vcenter{\xymatrix@-.5pc{
      \ar[r]|-@{|}^-{\fhat} \ar[d]_f \ar@{}[dr]|\Downarrow &
      \ar[r]|-@{|}^-{\ghat} \ar@{=}[d] \ar@{}[dr]|{1_{\ghat}} &
      \ar@{=}[d]\\
      \ar[r]|-{U_B} \ar[d]_g \ar@{}[dr]|{U_g} &
      \ar[r]|-{\ghat} \ar[d]|g \ar@{}[dr]|\Downarrow &
      \ar@{=}[d]\\
      \ar[r]|-@{|}_-{U_C} &
      \ar[r]|-@{|}_-{U_C} &
    }}\quad\text{and}\quad
  \vcenter{\xymatrix@-.5pc{
      \ar[r]|-@{|}^-{U_A} \ar@{=}[d] \ar@{}[dr]|\Downarrow &
      \ar[r]|-@{|}^-{U_A} \ar[d]|f \ar@{}[dr]|{U_f} &
      \ar[d]^f\\
      \ar[r]|-{\fhat} \ar@{=}[d] \ar@{}[dr]|{1_{\fhat}} &
      \ar[r]|-{U_B} \ar@{=}[d] \ar@{}[dr]|\Downarrow &
      \ar[d]^g\\
      \ar[r]|-@{|}_-{\fhat} &
      \ar[r]|-@{|}_-{\ghat} &
    }}
  \]
  It is easy to verify that these satisfy~\eqref{eq:compeqn}, using
  the interchange law for $\odot$ and $\circ$ in a double category.
\end{proof}

\begin{lem}\label{thm:theta-compose-horiz}
  Suppose that $f\maps A\to B$ has companions $\fhat$ and $\fhat'$,
  and that $g\maps B\to C$ has companions $\ghat$ and $\ghat'$.  Then
  $\theta_{\ghat,\ghat'}\odot \theta_{\fhat,\fhat'}  =
    \theta_{\ghat\odot\fhat, \ghat'\odot\fhat'}$.
\end{lem}
\begin{proof}
  Using the interchange law for $\odot$ and $\circ$, we have:
  \begin{align}
    \theta_{\ghat\odot\fhat, \ghat'\odot\fhat'} &=\;
    \vcenter{\xymatrix@-.5pc{
        \ar[r]|-@{|}^-{U_A} \ar@{=}[d] \ar@{}[dr]|\Downarrow &
        \ar[r]|-@{|}^-{U_A} \ar[d]|f \ar@{}[dr]|{U_f} &
        \ar[r]|-@{|}^-{\fhat} \ar[d]|f \ar@{}[dr]|\Downarrow &
        \ar[r]|-@{|}^-{\ghat} \ar@{=}[d] \ar@{}[dr]|{1_{\fhat}} &
        \ar@{=}[d]\\
        \ar[r]|-{\fhat'} \ar@{=}[d] \ar@{}[dr]|{1_{\ghat}} &
        \ar[r]|-{U_B} \ar@{=}[d] \ar@{}[dr]|\Downarrow &
        \ar[r]|-{U_B} \ar[d]|g \ar@{}[dr]|{U_g} &
        \ar[r]|-{\ghat} \ar[d]|g \ar@{}[dr]|\Downarrow &
        \ar@{=}[d]\\
        \ar[r]|-@{|}_-{\fhat'} &
        \ar[r]|-@{|}_-{\ghat'} &
        \ar[r]|-@{|}_-{U_C} &
        \ar[r]|-@{|}_-{U_C} &
      }}
    \;=\;
    \vcenter{\xymatrix@-.5pc{
        \ar[r]|-@{|}^-{U_A} \ar@{=}[d] \ar@{}[dr]|\Downarrow &
        \ar[r]|-@{|}^-{\fhat} \ar[d]|f \ar@{}[dr]|\Downarrow &
        \ar[r]|-@{|}^-{\ghat} \ar@{=}[d] \ar@{}[dr]|{1_{\fhat}} &
        \ar@{=}[d]\\
        \ar[r]|-{\fhat'} \ar@{=}[d] \ar@{}[dr]|{1_{\ghat}} &
        \ar[r]|-{U_B} \ar@{=}[d] \ar@{}[dr]|\Downarrow &
        \ar[r]|-{\ghat} \ar[d]|g \ar@{}[dr]|\Downarrow &
        \ar@{=}[d]\\
        \ar[r]|-@{|}_-{\fhat'} &
        \ar[r]|-@{|}_-{\ghat'} &
        \ar[r]|-@{|}_-{U_C} &
      }}\\
    &=\;
    \vcenter{\xymatrix@-.5pc{
        \ar[r]|-@{|}^-{U_A} \ar@{=}[d] \ar@{}[dr]|\Downarrow &
        \ar[r]|-@{|}^-{\fhat} \ar[d]|f \ar@{}[dr]|\Downarrow &
        \ar[r]|-@{|}^-{U_B} \ar@{=}[d] \ar@{}[dr]|{1_{U_B}} &
        \ar[r]|-@{|}^-{\ghat} \ar@{=}[d] \ar@{}[dr]|{1_{\fhat}} &
        \ar@{=}[d]\\
        \ar[r]|-{\fhat'} \ar@{=}[d] \ar@{}[dr]|{1_{\ghat}} &
        \ar[r]|-{U_B} \ar@{=}[d] \ar@{}[dr]|{1_{U_B}} &
        \ar[r]|-{U_B} \ar@{=}[d] \ar@{}[dr]|\Downarrow &
        \ar[r]|-{\ghat} \ar[d]|g \ar@{}[dr]|\Downarrow &
        \ar@{=}[d]\\
        \ar[r]|-@{|}_-{\fhat'} &
        \ar[r]|-@{|}_-{U_B} &
        \ar[r]|-@{|}_-{\ghat'} &
        \ar[r]|-@{|}_-{U_C} &
      }}\;=\;
    \vcenter{\xymatrix@-.5pc{
        \ar[r]|-@{|}^-{U_A} \ar@{=}[d] \ar@{}[dr]|\Downarrow &
        \ar[r]|-@{|}^-{\fhat} \ar[d]|f \ar@{}[dr]|\Downarrow &
        \ar[r]|-@{|}^-{U_B} \ar@{=}[d] \ar@{}[dr]|\Downarrow &
        \ar[r]|-@{|}^-{\ghat} \ar[d]|g \ar@{}[dr]|\Downarrow& \ar@{=}[d]\\
        \ar[r]|-@{|}_-{\fhat'} &
        \ar[r]|-@{|}_-{U_B} &
        \ar[r]|-@{|}_-{\ghat'} &
        \ar[r]|-@{|}_-{U_C} &
      }}\\
    &=\;
    \theta_{\ghat,\ghat'}\odot \theta_{\fhat,\fhat'} 
  \end{align}
  as desired.
\end{proof}

\begin{lem}\label{thm:theta-unit}
  If $f\maps A\to B$ has a companion \fhat, then
  $\theta_{\fhat,\fhat\odot U_A}$ and $\theta_{\fhat,U_B\odot \fhat}$
  are equal to the unit constraints $\fhat \iso \fhat\odot U_A$ and
  $\fhat\iso U_B\odot \fhat$.
\end{lem}
\begin{proof}
  By definition, we have
  \[\theta_{\fhat,\fhat\odot U_A} =\;
  \vcenter{\xymatrix@-.5pc{
      \ar[r]|-@{|}^-{U_A} \ar@{=}[d] \ar@{}[dr]|{\Downarrow 1_{U_A}} &
      \ar[r]|-@{|}^-{U_A} \ar@{=}[d] \ar@{}[dr]|{1_{U_A}} &
      \ar@{=}[d] \ar[rr]|-@{|}^-{\fhat} \ar@{}[ddrr]|\Downarrow && \ar@{=}[dd]\\
      \ar[r]|-{U_A} \ar@{=}[d] \ar@{}[dr]|{1_{U_A}} &
      \ar[r]|-{U_A} \ar@{=}[d] \ar@{}[dr]|\Downarrow &
      \ar[d]^f\\
      \ar[r]|-@{|}_-{U_A} &
      \ar[r]|-@{|}_-{\fhat} & \ar[rr]|-@{|}^-{U_B} &&
    }}\;=\;
  \vcenter{\xymatrix{ \ar[r]|-@{|}^-{U_A} \ar@{=}[d]
      \ar@{}[dr]|{\Downarrow 1_{U_A}} &  \ar@{=}[d]\\
      \ar[r]|-@{|}_-{U_A} & }}
  \]
  which, bearing in mind our suppression of unit and associativity
  constraints, means that in actuality it is the unit constraint
  $\fhat \iso \fhat\odot U_A$.  The other case is dual.
\end{proof}

\begin{lem}\label{thm:comp-func}
  Let $F\maps \lD\to\lE$ be a functor between double categories and
  let $f\maps A\to B$ have a companion \fhat\ in \lD.  Then $F(\fhat)$
  is a companion of $F(f)$ in \lE.
\end{lem}
\begin{proof}
  We take the defining 2-morphisms to be
  \[\vcenter{\xymatrix@R=1.5pc@C=3pc{
      \ar[r]|-@{|}^-{F(\fhat)} \ar[d]_{F(f)}
      \ar@{}[dr]|{F(\Downarrow)} &  \ar@{=}[d]\\
      \ar[r]|-{F(U_B)} \ar@{=}[d] \ar@{}[dr]|\iso &  \ar@{=}[d]\\
      \ar[r]|-@{|}_-{U_{F(B)}} & }}
  \quad\text{and}\quad
  \vcenter{\xymatrix@R=1.5pc@C=3pc{
      \ar[r]|-@{|}^-{U_{FA}} \ar@{=}[d] \ar@{}[dr]|\iso & \ar@{=}[d]\\
      \ar[r]|-{F(U_{A})} \ar@{=}[d] \ar@{}[dr]|{F(\Downarrow)} & 
      \ar[d]^{F(f)}\\
      \ar[r]|-@{|}_-{F(\fhat)} & .}}\]
  The axioms~\eqref{eq:compeqn} follow directly from those for \fhat.
\end{proof}

\begin{lem}\label{thm:comp-ten}
  Suppose that \lD\ is a monoidal double category and that $f\maps
  A\to B$ and $g\maps C\to D$ have companions \fhat\ and \ghat\
  respectively.  Then $\fhat\ten\ghat$ is a companion of $f\ten g$.
\end{lem}
\begin{proof}
  This follows from \autoref{thm:comp-func}, since $\ten\maps
  \lD\times\lD\to\lD$ is a functor, and a companion in $\lD\times\lD$
  is simply a pair of companions in \lD.
\end{proof}

\begin{lem}\label{thm:theta-func}
  Suppose that $f\maps \lD\to\lE$ is a functor and that $f\maps A\to
  B$ has companions \fhat\ and $\fhat'$ in \lD.  Then
  $\theta_{F(\fhat),F(\fhat')} = F(\theta_{\fhat,\fhat'})$.
\end{lem}
\begin{proof}
  Using the axioms of a pseudo double functor and the definition of
  the 2-morphisms in \autoref{thm:comp-func}, we have
  \begin{equation}
    F(\theta_{\fhat,\fhat'})
    =\;
    \vcenter{\xymatrix@C=4.5pc{
        \ar[r]|-@{|}^-{F(\fhat)}
        \ar[d] \ar@{}[dr]|{F(\Downarrow\odot\Downarrow)} &  \ar[d]\\
        \ar[r]|-@{|}_-{F(\fhat')} &}}
    \;=\;
    \vcenter{\xymatrix@C=2pc{
        \ar[rr]|-@{|}^-{F(\fhat)}
        \ar@{=}[d] \ar@{}[drr]|\iso &&  \ar@{=}[d]\\
        \ar[r]|-@{|}^-{F(U_{A})} \ar@{=}[d]
        \ar@{}[dr]|{F(\Downarrow)} &
        \ar[r]|-@{|}^-{F(\fhat)} \ar[d]|{F(f)}
        \ar@{}[dr]|{F(\Downarrow)}
        & \ar@{=}[d]\\
        \ar[r]|-@{|}_-{F(\fhat')} \ar@{}[drr]|\iso\ar@{=}[d] &
        \ar[r]|-@{|}_-{U_{F(B)}} & \ar@{=}[d]\\
        \ar[rr]|-@{|}_-{F(\fhat')} && }}
    \;=\;
    \vcenter{\xymatrix@R=1.5pc@C=2.5pc{
        \ar[r]|-@{|}^-{U_{F(A)}} \ar@{=}[d] \ar@{}[dr]|\iso &
        \ar[r]|-@{|}^-{F(\fhat)} \ar@{=}[d] \ar@{}[dr]|=
        & \ar@{=}[d]\\
        \ar[r]|-{F(U_{A})} \ar@{=}[d] \ar@{}[dr]|{F(\Downarrow)} &
        \ar[r]|-{F(\fhat)} \ar[d]|{F(f)} \ar@{}[dr]|{F(\Downarrow)}
        & \ar@{=}[d]\\
        \ar[r]|-{F(\fhat')}  \ar@{=}[d] \ar@{}[dr]|= &
        \ar[r]|-{F(U_{B})} \ar@{}[dr]|\iso  \ar@{=}[d] & \ar@{=}[d]\\
        \ar[r]|-@{|}_-{F(\fhat')} &
        \ar[r]|-@{|}_-{U_{F(B)}} &}}
    \;=
    \theta_{F(\fhat),\,F(\fhat')}
  \end{equation}
  as desired.
\end{proof}

\begin{lem}\label{thm:theta-ten}
  Suppose that \lD\ is a monoidal double category, that $f\maps A\to
  B$ has companions \fhat\ and $\fhat'$, and that $g\maps C\to D$ has
  companions \ghat\ and $\ghat'$.  Then $\theta_{\fhat,\fhat'} \ten
  \theta_{\ghat,\ghat'} = \theta_{\fhat\ten \ghat, \fhat'\ten\ghat'}.$
\end{lem}
\begin{proof}
  This follows from \autoref{thm:theta-func} in the same way that
  \autoref{thm:comp-ten} follows from \autoref{thm:comp-func}.
\end{proof}

\begin{lem}\label{thm:comp-iso}
  If $f\maps A\to B$ is a vertical isomorphism with a companion \fhat,
  then \fhat\ is a conjoint of its inverse $f\inv$.
\end{lem}
\begin{proof}
  The composites
  \[\vcenter{\xymatrix@-.5pc{
      \ar[r]|-@{|}^{\fhat}\ar[d]_f \ar@{}[dr]|{\Downarrow} &
      \ar@{=}[d]\\
      \ar[r]|{U_B}\ar[d]_{f\inv} \ar@{}[dr]|{\Downarrow U_{f\inv}} &
      \ar[d]^{f\inv}\\
      \ar[r]|-@{|}_{U_A} &
    }}\quad\text{and}\quad
  \vcenter{\xymatrix@-.5pc{
      \ar[r]|-@{|}^{U_B}\ar[d]_{f\inv} \ar@{}[dr]|{\Downarrow U_{f\inv}} &
      \ar[d]^{f\inv}\\
      \ar[r]|{U_A}\ar@{=}[d] \ar@{}[dr]|{\Downarrow} &
      \ar[d]^f\\
      \ar[r]|-@{|}_{\fhat} &
    }}
  \]
  exhibit \fhat\ as a conjoint of $f\inv$.
\end{proof}

\begin{lem}\label{thm:compconj-adj}
  If $f\maps A\to B$ has both a companion \fhat\ and a conjoint \fchk,
  then we have an adjunction $\fhat\adj\fchk$ in $\cH\lD$.  If $f$ is
  an isomorphism, then this is an adjoint equivalence.
\end{lem}
\begin{proof}
  The unit and counit of the adjunction $\fhat\adj\fchk$ are the
  composites
  \[\vcenter{\xymatrix@-.5pc{
      \ar[r]|-@{|}^{U_A}\ar@{=}[d] \ar@{}[dr]|{\Downarrow} &
      \ar[r]|-@{|}^{U_A}\ar[d]|{f} \ar@{}[dr]|{\Downarrow} &
      \ar@{=}[d]\\
      \ar[r]|-@{|}_{\fhat} &
      \ar[r]|-@{|}_{\fchk} &
    }}\quad\text{and}\quad
  \vcenter{\xymatrix@-.5pc{
      \ar[r]|-@{|}^{\fchk}\ar@{=}[d] \ar@{}[dr]|{\Downarrow} &
      \ar[r]|-@{|}^{\fhat}\ar[d]|{f} \ar@{}[dr]|{\Downarrow} &
      \ar@{=}[d]\\
      \ar[r]|-@{|}_{U_B} &
      \ar[r]|-@{|}_{U_B} &
    }}
  \]
  The triangle identities follow from~\eqref{eq:compeqn}.  If $f$ is
  an isomorphism, then by the dual of \autoref{thm:comp-iso}, \fchk\
  is a companion of $f\inv$.  But then by \autoref{thm:comp-compose}
  $\fchk\odot \fhat$ is a companion of $1_A=f\inv \circ f$ and
  $\fhat\odot\fchk$ is a companion of $1_B = f\circ f\inv$, and hence
  \fhat\ and \fchk\ are equivalences.  We can then check that in this
  case the above unit and counit actually are the isomorphisms
  $\theta$, or appeal to the general fact that any adjunction
  involving an equivalence is an adjoint equivalence.
\end{proof}

\begin{rmk}
  Our intended applications actually only require our double
  categories to have companions and conjoints for vertical
  \emph{isomorphisms}; we may call a double category with this
  property \textbf{isofibrant}.  Note that by \autoref{thm:comp-iso},
  having companions for all isomorphisms implies having conjoints for
  all isomorphisms.  However, most examples we are interested in have
  all companions and conjoints, and these are useful for other
  purposes as well; see~\cite{shulman:frbi}.  Moreover, if we are
  given a double category in which only vertical isomorphisms have
  companions, we can still apply our theorems to it as written, simply
  by first discarding all noninvertible vertical 1-morphisms.
\end{rmk}

\section{From double categories to bicategories}
\label{sec:1x1-to-bicat}

We are now equipped to lift structures on fibrant double categories to
their horizontal bicategories.  In this section we show that passage
from fibrant double categories to bicategories is functorial; in the
next section we show that it preserves monoidal structure.

As a point of notation, we write $\odot$ for the composition of
1-cells in a bicategory, since our bicategories are generally of the
form $\cH(\lD)$.  As advocated by Max Kelly, we say \textbf{functor}
to mean a morphism between bicategories that preserves composition up
to isomorphism; equivalent terms include \emph{weak 2-functor},
\emph{pseudofunctor}, and \emph{homomorphism}.

\begin{thm}
  If \lD\ is a double category, then $\cH(\lD)$ is a bicategory, and
  any functor $F\maps \lD\to\lE$ induces a functor $\cH(F)\maps
  \cH(\lD)\to\cH(\lE)$.  In this way $\cH$ defines a functor of
  1-categories $\mathbf{Dbl}\to \mathbf{Bicat}$.
\end{thm}
\begin{proof}
  The constraints of $F$ are all globular, hence give constraints for
  $\cH(F)$.  Functoriality is evident.
\end{proof}

The action of \cH\ on transformations, however, is less obvious, and
requires the presence of companions or conjoints.  Recall that if
$F,G\maps \cA\to\cB$ are functors between bicategories, then an
\textbf{oplax transformation} $\al\maps F\to G$ consists of 1-cells
$\al_A\maps FA\to GA$ and 2-cells
\[\vcenter{\xymatrix{ \ar[r]^{Ff}\ar[d]_{\al_A} \drtwocell\omit{\al_f} &  \ar[d]^{\al_B}\\
  \ar[r]_{Gf} & }}\]
such that for any 2-cell $\xymatrix{A \rtwocell^f_g{x} & B}$ in \cA,
\begin{equation}
  \label{eq:laxtransf-nat}
  \vcenter{\xymatrix@R=1pc@C=3pc{
      \rtwocell^{Ff}_{Fg}{Fx}\ar[dd]_{\al_A} 
      &  \ar[dd]^{\al_B}\\
      \drtwocell\omit{\al_g} & \\
      \ar[r]_{Gg} & }}\;=\;
  \vcenter{\xymatrix@R=1pc@C=3pc{
      \ar[r]^{Ff}\ar[dd]_{\al_A} \drtwocell\omit{\al_f} &
      \ar[dd]^{\al_B}\\ & \\
      \rtwocell^{Gf}_{Gg}{Gx} & }}
\end{equation}
and moreover for any $A$ and any $f,g$ in \cA,
\begin{equation}
  \vcenter{\xymatrix@R=5pc{
      \rtwocell^{1_{FA}}_{F(1_A)}{\iso} \ar[d]_{\al_A} \drtwocell\omit{\al_{1_A}} &  \ar[d]^{\al_A}\\
      \rtwocell^{G(1_A)}_{1_{GA}}{\iso} & }} \;=\;
  \vcenter{\xymatrix{ \ar[r]^{1_{FA}}\ar[d]_{\al_A} \drtwocell\omit{\iso}&  \ar[d]^{\al_A}\\
      \ar[r]_{1_{GA}} &
    }}
  \quad\text{and}\quad
  \vcenter{\xymatrix{
      \ar[r]|{Ff}\ar[d]_{\al_A} \drtwocell\omit{\al_f}
      \rruppertwocell^{F(gf)}{\iso}
      &
      \ar[r]|{Fg}\ar[d]|{\al_B} \drtwocell\omit{\al_g} &
      \ar[d]^{\al_C}\\
      \ar[r]|{Gf} \rrlowertwocell_{G(gf)}{\iso} & \ar[r]|{Gg} & }}
  \;=\;
  \vcenter{\xymatrix{ \ar[r]^{F(gf)}\ar[d]_{\al_A} \drtwocell\omit{\al_{gf}} &  \ar[d]^{\al_C}\\
      \ar[r]_{G(gf)} & }}\label{eq:laxtransf-ax}
\end{equation}
It is a \textbf{lax transformation} if the 2-cells $\al_f$ go the
other direction, and a \textbf{pseudo transformation} if they are
isomorphisms.

By doctrinal adjunction~\cite{kelly:doc-adjn}, given collections of
1-cells $\al_A\maps FA\to GA$ and $\be_A\maps GA\to FA$ and
adjunctions $\al_A\adj \be_A$ in \cB, there is a bijection between
\begin{inparaenum}
\item collections of 2-cells $\al_f$ making $\al$ an oplax
  transformation and
\item collections of 2-cells $\be_f$ making $\be$ a lax
  transformation.
\end{inparaenum}
Two such transformations correspond under this bijection if and only if
\begin{equation}
  \vcenter{\xymatrix@-.5pc{F(f) \ar[r]^-{\eta \odot F(f)}
      \ar[d]_{F(f)\odot \eta} &
      \be_B\odot \al_B \odot F(f) \ar[d]^{\be_B \odot \al_f}\\
      F(f) \odot \be_A\odot \al_A\ar[r]_-{\be_f \odot \al_A} &
      \be_B\odot G(f) \odot \al_A}}
  \quad\text{and}\quad
  \vcenter{\xymatrix@-.5pc{\al_B\odot F(f)\odot \be_A
      \ar[r]^-{\al_B\odot \be_f}\ar[d]_{\al_f \odot \be_A}&
      \al_B \odot \be_B \odot G(f)\ar[d]^{\ep \odot G(f)}\\
      G(f)\odot \al_A\odot \be_A \ar[r]_-{G(f) \odot \ep} & G(f)}}\label{eq:conjtrans}
\end{equation}
commute.  If we have a pointwise adjunction between an oplax and a lax
transformation, whose 2-cell structures correspond under this
bijection, we call it a \textbf{conjunctional transformation}
$(\al\conj \be)\maps F\to G$.  (These are the conjoint pairs in a
double category whose horizontal arrows are lax transformations and
whose vertical arrows are oplax transformations.)

Of particular importance is the case when both $\al$ and \be\ are
pseudo natural and each adjunction $\al_A\adj \be_A$ is an adjoint
equivalence.  In this case we call $\al\conj \be$ a \textbf{pseudo
  natural adjoint equivalence}.  A pseudo natural adjoint equivalence
can equivalently be defined as an internal equivalence in the
bicategory $\cBicat(\cA,\cB)$ of functors, pseudo natural
transformations, and modifications $\cA\to\cB$.

Recall also that if $\al,\al'\maps F\to G$ are oplax transformations,
a \textbf{modification} $\mu\maps \al\to\al'$ consists of 2-cells
$\mu_A\maps \al_A\to\al'_A$ such that
\begin{equation}
  \vcenter{\xymatrix@C=1pc@R=2.5pc{ \ar[rr]^{Ff}\dtwocell_{\al'_A}^{\al_A}{\mu_A}  &
      \drtwocell\omit{\al_f} &  \ar[d]^{\al_B}\\
      \ar[rr]_{Gf} && }} \quad=\quad
  \vcenter{\xymatrix@C=1pc@R=2.5pc{ \ar[rr]^{Ff}\ar[d]_{\al'_A} \drtwocell\omit{\al'_f} && 
      \dtwocell^{\al_B}_{\al'_B}{\mu_B}\\
      \ar[rr]_{Gf} && }}\label{eq:modif-ax}
\end{equation}
There is an evident notion of modification between lax transformations
as well.  Finally, given conjunctional transformations $\al\conj\be$
and $\al'\conj \be'$, there is a bijection between modifications
$\al\to\al'$ and $\be'\to\be$, where $\mu\maps \al\to\al'$ corresponds
to $\bar{\mu}\maps \be'\to\be$ with components $\bar{\mu}_A$ defined
by:
\[\vcenter{\xymatrix@-.5pc{
    && FA \ar@{=}[drr] \ddtwocell<5>^{\al_A}_{\al'_A}{\mu_A}\\
    GA \ar[urr]^{\be'_A} \ar@{=}[drr] & \Swarrow_\ep && \Swarrow_\eta & FA\\
    &&GA\ar[urr]_{\be_A}
  }}\]
The modifications $\bar{\mu}$ and \mu\ are called \textbf{mates}, and
are compatible with composition (see \cite{ks:r2cats}).  Thus, given
$\cA,\cB$ we can define a bicategory $\Conj(\cA,\cB)$, whose objects
are functors $\cA\to\cB$, whose 1-cells are conjunctional
transformations considered as pointing in the direction of their left
adjoints, and whose 2-cells are mate-pairs of modifications.

\begin{thm}\label{thm:h-locfr}
  If \lD\ is a double category and \lE\ is a fibrant double category
  with chosen companions and conjoints, we have a functor
  \begin{align}
    \cDbl(\lD,\lE) &\too \Conj(\cH(\lD),\cH(\lE))\\
    F &\mapsto \cH(F)\\
    \al &\mapsto (\alhat\conj\alchk).
  \end{align}
  Moreover, if \al\ is an isomorphism, then $\alhat\conj\alchk$ is a
  pseudo natural adjoint equivalence.
\end{thm}

Note that we are here regarding the 1-category $\cDbl(\lD,\lE)$ as a
bicategory with only identity 2-cells.

\begin{proof}
  We denote the chosen companion and conjoint of $f$ in \lE\ by \fhat\
  and \fchk, as usual.  We define $\alhat$ as follows: its 1-cell
  components are $\alhat_A = \widehat{\al_A}$, and its 2-cell
  component $\alhat_f$ is the composite
  \begin{equation}
    \vcenter{\xymatrix@R=1.5pc@C=2.5pc{
        \ar[r]|-@{|}^{U_{FA}}\ar@{=}[d] \ar@{}[dr]|{\Downarrow} &
        \ar[r]^{Ff}\ar[d]|{\al_A} \ar@{}[dr]|{\Downarrow \al_f} &
        \ar[r]|-@{|}^{\alhat_B}\ar[d]|{\al_B} \ar@{}[dr]|{\Downarrow} &
        \ar@{=}[d]\\
        \ar[r]|-@{|}_{\alhat_A} &
        \ar[r]_{Gf} &
        \ar[r]|-@{|}_{U_{GB}} & 
      }}\label{eq:oplax-2cell}
  \end{equation}
  Equations~\eqref{eq:laxtransf-nat} and~\eqref{eq:laxtransf-ax}
  follow directly from \autoref{thm:dbl-transf}.  The construction of
  $\alchk$ is dual, using conjoints, and \autoref{thm:compconj-adj}
  shows that $\alhat_A\adj \alchk_A$.  For the first equation
  in~\eqref{eq:conjtrans}, we have
  \begin{equation}
    \vcenter{\xymatrix@-.5pc{
        \ar[r]|-@{|}^{U_{FA}}\ar@{=}[d] \ar@{}[dr]|{=} &
        \ar[r]|-@{|}^{Ff}\ar@{=}[d] \ar@{}[dr]|{=} &
        \ar[r]|-@{|}^{U_{FB}}\ar@{=}[d] \ar@{}[dr]|{\Downarrow} &
        \ar[r]|-@{|}^{U_{FB}}\ar[d]|{\al_B} \ar@{}[dr]|{\Downarrow} &
        \ar@{=}[d]\\
        \ar[r]|{U_{FA}}\ar@{=}[d] \ar@{}[dr]|{\Downarrow} &
        \ar[r]|{Ff}\ar[d]|{\al_A} \ar@{}[dr]|{\Downarrow\al_f} &
        \ar[r]|{\alhat_B}\ar[d]|{\al_B} \ar@{}[dr]|{\Downarrow} &
        \ar[r]|{\alchk_B}\ar@{=}[d] \ar@{}[dr]|{=} &
        \ar@{=}[d]\\
        \ar[r]|-@{|}_{\alhat_A} &
        \ar[r]|-@{|}_{Gf} &
        \ar[r]|-@{|}_{U_{GB}} &
        \ar[r]|-@{|}_{\alchk_B} &
      }}\;=\;
    \vcenter{\xymatrix@-.5pc{
        \ar[r]|-@{|}^{U_{FA}}\ar@{=}[d] \ar@{}[dr]|{\Downarrow} &
        \ar[r]|-@{|}^{Ff}\ar[d]|{\al_A} \ar@{}[dr]|{\Downarrow \al_f} &
        \ar[r]|-@{|}^{U_{FB}}\ar[d]|{\al_B} \ar@{}[dr]|{\Downarrow} &
        \ar@{=}[d]\\
        \ar[r]|-@{|}_{\alhat_A} &
        \ar[r]|-@{|}_{Gf} &
        \ar[r]|-@{|}_{\alchk_B} &
      }}\;=\;
    \vcenter{\xymatrix@-.5pc{
        \ar[r]|-@{|}^{U_{FA}}\ar@{=}[d] \ar@{}[dr]|{\Downarrow} &
        \ar[r]|-@{|}^{U_{FA}}\ar[d]|{\al_A} \ar@{}[dr]|{\Downarrow} &
        \ar[r]|-@{|}^{Ff}\ar@{=}[d] \ar@{}[dr]|{=} &
        \ar[r]|-@{|}^{U_{FB}}\ar@{=}[d] \ar@{}[dr]|{=} &
        \ar@{=}[d]\\
        \ar[r]|{\alhat_A}\ar@{=}[d] \ar@{}[dr]|{=} &
        \ar[r]|{\alchk_A}\ar@{=}[d] \ar@{}[dr]|{\Downarrow} &
        \ar[r]|{Ff}\ar[d]|{\al_A} \ar@{}[dr]|{\Downarrow\al_f} &
        \ar[r]|{U_{FB}}\ar[d]|{\al_B} \ar@{}[dr]|{\Downarrow} &
        \ar@{=}[d]\\
        \ar[r]|-@{|}_{\alhat_A} &
        \ar[r]|-@{|}_{U_{GA}} &
        \ar[r]|-@{|}_{Gf} &
        \ar[r]|-@{|}_{\alchk_B} &
        }},
  \end{equation}
  and the second is dual.  Thus $(\alhat\conj\alchk)$ is a
  conjunctional transformation.

  Now suppose given $\al\maps F\to G$ and $\be\maps G\to H$.  Then by
  \autoref{thm:comp-compose}, $\behat_A\odot\alhat_A$ is a companion
  of $\be_A\circ \al_A$, so we have a canonical isomorphism
  \[\theta_{\widehat{\be\al}_A, \,\behat_A\odot\alhat_A}\maps
  \widehat{\be\al}_A \too[\iso] \behat_A\odot\alhat_A.
  \]
  Of course, we also have $\theta_{\widehat{1_A},U_A}\maps
  \widehat{1_A} \too[\iso] U_A$ by \autoref{thm:comp-unit}.  These
  constraints are automatically natural, since $\cDbl(\lD,\lE)$ has no
  nonidentity 2-cells.  The axiom for the composition constraint says
  that two constructed isomorphisms
  \[\widehat{\gm\be\al}_A \too[\iso] (\gmhat_A \odot \behat_A)\odot \alhat_A\]
  are equal.  However, both $\widehat{\gm\be\al}_A$ and $(\gmhat_A
  \odot \behat_A)\odot \alhat_A$ are companions of $\gm_A\be_A\al_A$,
  and both of these isomorphisms are constructed from composites (both
  $\circ$-composites and $\odot$-composites) of $\theta$s; hence by
  Lemmas \ref{thm:theta-compose-vert} and
  \ref{thm:theta-compose-horiz} they are both equal to
  \[\theta_{\widehat{\gm\be\al}_A,\, (\gmhat_A \odot \behat_A)\odot
    \alhat_A}\] and thus equal to each other.  The same argument
  applies to the axioms for the unit constraint; thus we have a
  functor of bicategories.

  Finally, if $\al$ is an isomorphism, then in particular each $\al_A$
  is an isomorphism, so by \autoref{thm:compconj-adj} each
  $\alhat_A\adj \alchk_A$ is an adjoint equivalence.  But \al\ being
  an isomorphism also implies that each 2-cell
  \[\vcenter{\xymatrix@-.5pc{ \ar[r]|-@{|}^-{Ff} \ar[d]_{\al_A} \ar@{}[dr]|{\Downarrow\al_f} &  \ar[d]^{\al_B}\\
      \ar[r]|-@{|}_-{Gf} & }}\]
  is an isomorphism.  From its inverse we form the composite
  \[\vcenter{\xymatrix@R=1.5pc@C=3pc{
      \ar[r]|-@{|}^{\alhat_A}\ar@{=}[d] \ar@{}[dr]|{\Downarrow} &
      \ar[r]^{Gf}\ar[d]|{\al_A\inv} \ar@{}[dr]|{\Downarrow\al_f\inv} &
      \ar[r]|-@{|}^{U_{GB}}\ar[d]|{\al_B\inv} \ar@{}[dr]|{\Downarrow} &
      \ar@{=}[d]\\
      \ar[r]|-@{|}_{U_{FA}}&
      \ar[r]_{Ff} &
      \ar[r]|-@{|}_{\alchk_A} &
    }}\]
  which we can then verify to be an inverse of~\eqref{eq:oplax-2cell}.
  Thus $\alhat$, and dually $\alchk$, is pseudo natural, and hence
  $\alhat\conj\alchk$ is a pseudo natural adjoint equivalence.
\end{proof}

We can also promote \autoref{thm:theta} to a functorial uniqueness.

\begin{lem}\label{thm:h-locfr-uniq}
  Let \lD\ be a double category and \lE\ a fibrant double category
  with two different sets of choices $\fhat,\fchk$ and $\fhat',\fchk'$
  of companions and conjoints for each vertical 1-morphism $f$, giving
  rise to two different functors
  \[\cH,\cH'\maps \cDbl(\lD,\lE)\too \Conj(\cH(\lD),\cH(\lE)).\]
  Then the isomorphisms $\theta$ from \autoref{thm:theta} fit together
  into a pseudo natural adjoint equivalence $\cH\eqv \cH'$ which is the
  identity on objects.
\end{lem}
\begin{proof}
  We must first show that for a given transformation $\al\maps F\to
  G\maps \lD\to\lE$ in \cDbl, the isomorphisms \th\ form an invertible
  modification $\alhat \iso \alhat'$.
  Substituting~\eqref{eq:oplax-2cell} and the definition of \th\
  into~\eqref{eq:modif-ax}, this becomes the assertion that
  \begin{equation}
    \vcenter{\xymatrix@R=1.5pc@C=2pc{
        &
        \ar[r]|-@{|}^{U_{FA}}\ar@{=}[d] \ar@{}[dr]|{\Downarrow} &
        \ar[r]^{Ff}\ar[d]|{\al_A} \ar@{}[dr]|{\Downarrow \al_f} &
        \ar[r]|-@{|}^{\alhat_B}\ar[d]|{\al_B} \ar@{}[dr]|{\Downarrow} &
        \ar@{=}[d]\\
        \ar[r]|-@{|}^{U_{FA}} \ar@{=}[d] \ar@{}[dr]|{\Downarrow} &
        \ar[r]|{\alhat_A} \ar[d]|{\al_A} \ar@{}[dr]|{\Downarrow}&
        \ar[r]_{Gf}  \ar@{=}[d] &
        \ar[r]|-@{|}_{U_{GB}} & \\
        \ar[r]|-@{|}_{\alhat_A'} & \ar[r]|-@{|}_{U_{GB}}&&
      }} \;=\;
    \vcenter{\xymatrix@R=1.5pc@C=2pc{
        && \ar@{=}[d] \ar[r]|-@{|}^{U_{FA}} \ar@{}[dr]|{\Downarrow} &
        \ar[d]|{\al_B} \ar[r]|-@{|}^{\alhat_B} \ar@{}[dr]|{\Downarrow}
        &
        \ar@{=}[d] &\\
        \ar[r]|-@{|}^{U_{FA}}\ar@{=}[d] \ar@{}[dr]|{\Downarrow} &
        \ar[r]^{Ff}\ar[d]|{\al_A} \ar@{}[dr]|{\Downarrow \al_f} &
        \ar[r]|{\alhat_B'}\ar[d]|{\al_B} \ar@{}[dr]|{\Downarrow} &
        \ar@{=}[d] \ar[r]|-@{|}_{U_{GB}}&\\
        \ar[r]|-@{|}_{\alhat_A'} &
        \ar[r]_{Gf} &
        \ar[r]|-@{|}_{U_{GB}} & .
      }}
  \end{equation}
  This follows from two applications of~\eqref{eq:compeqn}, one for
  $\alhat_A$ and one for $\alhat_B'$.  (The mate of \th\ is, of
  course, uniquely determined.)  Now, to show that these form a pseudo
  natural adjoint equivalence, it remains only to check that they do,
  in fact, form a pseudo natural transformation which is the identity
  on objects, i.e.\ that~\eqref{eq:laxtransf-nat}
  and~\eqref{eq:laxtransf-ax} are satisfied.
  But~\eqref{eq:laxtransf-nat} is vacuous since $\cDbl(\lD,\lE)$ has
  no nonidentity 2-cells, and~\eqref{eq:laxtransf-ax} follows from
  Lemmas \ref{thm:theta-compose-vert} and
  \ref{thm:theta-compose-horiz} since all the constraints involved are
  also instances of \th.
\end{proof}

It seems that we should have a functor from fibrant double categories
to a tricategory of bicategories, functors, conjunctional
transformations, and modifications, but there is no tricategory
containing conjunctional transformations since the interchange law
only holds laxly.  However, we can say the following.  Let $\cDbl^f_g$
denote the sub-2-category of \cDbl\ containing the fibrant double
categories, all functors between them, and only the transformations
that are isomorphisms, and let \cBicat\ denote the tricategory of
bicategories, functors, pseudo natural transformations, and
modifications.

\begin{thm}\label{thm:h-functor}
  There is a functor of tricategories $\cH\maps \cDbl^f_g\to \cBicat$.
\end{thm}
\begin{proof}
  The definition of functors between tricategories can be found
  in~\cite{gps:tricats} or~\cite{nick:tricats}.  In addition to
  \autoref{thm:h-locfr}, we require pseudo natural (adjoint)
  equivalences $\chi$ and $\iota$ relating composition and units in
  $\cDbl^f_g$ and \cBicat, and modifications relating composites of
  these, which satisfy various axioms.  However, since composition of
  1-cells in $\cDbl^f_g$ and \cBicat\ is strictly associative and
  unital, \cH\ strictly preserves this composition, and $\cDbl^f_g$
  has no nonidentity 3-cells, this merely amounts to the following.

  Firstly, for every pair of transformations
  \[\vcenter{\xymatrix{\lC \rtwocell^F_G{\al} & \lD \rtwocell^H_K{\be}
      & \lE}}\]
  between fibrant double categories, we require an invertible
  modification $\chi\maps \behat * \alhat \iso \widehat{\be*\al}$ such
  that
  \begin{equation}
    \vcenter{\xymatrix@-.5pc{1 \ar[r] \ar[dr] & \hat{1}*\hat{1} \ar[d]^\chi\\
        & \widehat{1*1} }} \quad\text{and}\quad
    \vcenter{\xymatrix@-.5pc{
        \widehat{\gm\al}*\widehat{\de\be} \ar[r]\ar[d]_\chi &
        (\gmhat*\dehat)(\alhat*\behat)\ar[d]^{\chi\chi}\\
        \widehat{\gm\al*\de\be}\ar[r] &
        (\widehat{\gm*\de})(\widehat{\al*\be})}}
  \end{equation}
  commute.  (Here we are writing $*$ for the `Godement product' of
  2-cells in $\cDbl$ and $\cBicat$.)  These are the 2-cell
  components of the composition constraint, its 1-cell components
  being identities.  Now by Lemmas \ref{thm:comp-compose} and
  \ref{thm:comp-func}, $(\behat *\alhat)_A = \behat_{GA} \circ
  H(\alhat_A)$ is a companion of $(\be*\al)_A = \be_{GA} \circ
  H(\al_A)$.  Therefore, we take the component $\chi_A$ to be
  \[\theta_{\behat_{GA} \circ H(\alhat_A),\, \widehat{\be*\al}_A}.\]
  Equation~\eqref{eq:modif-ax}, saying that these form a modification,
  becomes the equality of two large composites of 2-cells in \lD,
  which as usual follows from~\eqref{eq:compeqn}.

  Secondly, for every $F\maps \lD\to\lE$ we require an isomorphism
  $\iota\maps \widehat{1_F} \iso 1_{\cH(F)}$ satisfying a couple of
  axioms which simply require it to be equal to the unit constraint of
  the local functor \cH\ from \autoref{thm:h-locfr}; these are the
  2-cell components of the unit constraint.  Finally, the required
  modifications merely amount to the \emph{assertions} that
  \[\vcenter{\xymatrix@-.5pc{\gmhat*\behat*\alhat \ar[r]^\chi\ar[d]_\chi &
      \widehat{\gm*\be}*\alhat \ar[d]^\chi\\
      \gmhat*\widehat{\be*\al}\ar[r]_\chi & \widehat{\gm*\be*\al}}},\qquad
  \vcenter{\xymatrix@-.5pc{ \alhat \ar[r]^-\iota \ar@{=}[dr] &
      \widehat{1_F}*\alhat \ar[d]^\chi \\ & \alhat }}, \;\text{and}\qquad
  \vcenter{\xymatrix@-.5pc{ \alhat \ar[r]^-\iota \ar@{=}[dr] &
      \alhat*\widehat{1_F} \ar[d]^\chi \\ & \alhat }}\]
  commute; again this follows from \autoref{thm:theta-compose-vert}.
\end{proof}

We end this section with one final lemma.

\begin{lem}\label{thm:theta-nat}
  Suppose $F,G\maps \lD\to\lE$ are functors, $\al\maps F\to G$ is a
  transformation, and that $f\maps A\to B$ has a companion \fhat\ in
  \lD.  Then the oplax comparison 2-cell for \alhat:
  \[\vcenter{\xymatrix{
      \ar[r]^{F(\fhat)}\ar[d]_{\alhat_A} \drtwocell\omit{\;\alhat_{\fhat}}&  \ar[d]^{\alhat_B}\\
      \ar[r]_{G(\fhat)} & }}\]
  is equal to $\theta_{\alhat_B\odot F(\fhat),\, G(\fhat) \odot
    \alhat_A}$ (and in particular is an isomorphism).
\end{lem}
\begin{proof}
  By definition $\alhat_A$ and $\alhat_B$ are companions of $\al_A$
  and $\al_B$, respectively, and by \autoref{thm:comp-func} $F(\fhat)$
  and $G(\fhat)$ are companions of $F(f)$ and $G(f)$, respectively.
  Thus, by \autoref{thm:comp-compose} the domain and codomain of
  $\alhat_{\fhat}$ are both companions of $G(f) \circ \al_A = \al_B
  \circ F(f)$, so at least the asserted $\theta$ isomorphism exists.
  Now, by taking the definition~\eqref{eq:oplax-2cell} of
  $\alhat_{\fhat}$ and substituting it for $\theta$
  in~\eqref{eq:comp-iso}, using the axioms for companions and the
  naturality of $\al$ on 2-morphisms, we see that $\alhat_{\fhat}$
  satisfies~\eqref{eq:comp-iso} and hence must be equal to $\theta$.
\end{proof}

\section{Symmetric monoidal bicategories}
\label{sec:constr-symm-mono}

We are now ready to lift monoidal structures from double categories to
bicategories.  If we had a theory of symmetric monoidal tricategories,
we could do this by improving \autoref{thm:h-functor} to say that
$\cH$ is a symmetric monoidal functor, and then conclude that it
preserves pseudomonoids.  However, in the absence of such a theory, we
give a direct proof.

\begin{thm}\label{thm:mon11-monbi}
  If \lD\ is a fibrant monoidal double category, then $\cH(\lD)$ is a
  monoidal bicategory.  If \lD\ is braided, so is $\cH(\lD)$, and if
  \lD\ is symmetric, so is $\cH(\lD)$.
\end{thm}

\begin{rmk}
  For monoidal bicategories, there is a notion in between braided and
  symmetric, called \emph{sylleptic}, in which the the braiding is
  self-inverse up to an isomorphism (the \emph{syllepsis}) but this
  isomorphism is not maximally coherent.  Since in our approach the
  syllepsis will be an isomorphism of the form
  $\theta_{\fhat,\fhat'}$, it is \emph{always} maximally coherent;
  thus our method cannot produce sylleptic monoidal bicategories that
  are not symmetric.
\end{rmk}

\begin{proof}[Proof of \autoref{thm:mon11-monbi}]
  A monoidal bicategory is defined to be a tricategory with one
  object.  We use the definition of tricategory
  from~\cite{nick:tricats}, which is the same as that
  of~\cite{gps:tricats} except that the associativity and unit
  constraints are pseudo natural adjoint equivalences, rather than
  merely pseudo transformations whose components are equivalences.

  The functor \cH\ evidently preserves products, so $\ten\maps
  \lD\times\lD\to\lD$ induces a functor $\ten\maps
  \cH(\lD)\times\cH(\lD)\to \cH(\lD)$, and of course $I$ is still the
  unit.  The associativity constraint of \lD\ is a natural isomorphism
  \[\vcenter{\xymatrix@C=5pc{\lD\times\lD\times\lD \rtwocell^{\ten
        (\Id\times\ten)}_{\ten(\ten\times\Id)}{\fa\iso} &\lD }}\]
  so by \autoref{thm:h-locfr} it gives rise to a pseudo natural
  adjoint equivalence
  \[\vcenter{\xymatrix@C=6pc{\cH(\lD)\times\cH(\lD)\times\cH(\lD) \rtwocell^{\ten
        (\Id\times\ten)}_{\ten(\ten\times\Id)}{\fahat\eqv} &\cH(\lD) }}\]
  Likewise, the unit constraints of \lD\ induce pseudo natural adjoint
  equivalences.

  The final four pieces of data for a monoidal bicategory are
  invertible modifications relating various composites of the
  associativity and unit transformations.  The first is a
  ``pentagonator'' which relates the two ways to go around the Mac
  Lane pentagon:
  \[\xy
  (-10,0)*{((A\ten B)\ten C)\ten D}="A";
  (20,10)*{(A\ten (B\ten C))\ten D}="B";
  (50,0)*{A\ten ((B\ten C)\ten D)}="C";
  (0,-15)*{(A\ten B)\ten (C\ten D)}="D";
  (40,-15)*{A\ten (B\ten (C\ten D))}="E";
  (20,-5)*{\scriptstyle\pi\Downarrow\iso};
  \ar "B";"A";^{\fahat\ten U_D}
  \ar "C";"B";^{\fahat}
  \ar "D";"A";_{\fahat}
  \ar "E";"D";_{\fahat}
  \ar "E";"C";^{U_A\ten \fahat}
  \endxy
  \]
  Now by Lemmas \ref{thm:comp-compose} and \ref{thm:comp-ten}, both
  sides of this pentagon in $\cH(\lD)$ are companions of the
  corresponding sides of the pentagon in $\lD_0$.  Since the pentagon
  in $\lD_0$ commutes, we have an isomorphism $\theta$ between the two
  sides of the pentagon in $\cH(\lD)$, which we take to be $\pi$.
  That \pi\ is in fact a modification follows from
  \autoref{thm:h-locfr-uniq}.  We construct the other invertible
  modifications $\mu, \lambda, \rho$ in the same way.

  Finally, we must show that three equations between pasting
  composites of 2-cells hold, relating composites of
  $\pi,\mu,\lambda,\rho$.  However, in each of these equations, both
  the domain and the codomain of the 2-cells involved are companions
  of the same isomorphism in $\lD_0$.  For the 5-associahedron, this
  isomorphism is the unique constraint
  \[(((A\ten B)\ten C)\ten D)\ten E \too[\iso] A\ten (B\ten (C\ten
  (D\ten E)));
  \]
  for the other two it is simply the associator $(A\ten B)\ten C
  \too[\iso] A\ten (B\ten C)$.  By Lemmas \ref{thm:theta-unit},
  \ref{thm:theta-ten}, and \ref{thm:theta-nat},
  every 2-cell in these diagrams is a $\theta$ isomorphism relating
  two companions of the same vertical isomorphism.  Therefore, Lemmas
  \ref{thm:theta-compose-vert} and \ref{thm:theta-compose-horiz} imply
  that each pasting diagram is also a $\theta$ isomorphism between its
  domain and codomain.  The uniqueness of $\theta$ then implies that
  the three equations hold.

  Now suppose that \lD\ is braided; to show that $\cH(\lD)$ is braided
  we seemingly must first have a definition of braided monoidal
  bicategory.  The interested reader may follow the tortuous path of
  the definition of braided monoidal 2-categories and bicategories
  through the literature, starting from~\cite{kv:2cat-zam,kv:bm2cat}
  and continuing, with occasional corrections,
  through~\cite{bn:hda-i,ds:monbi-hopfagbd,crans:centers,mccrudden:bal-coalgb},
  and~\cite{gurski:brmonbicat}.  However, the details of the
  definition are essentially unimportant for us; since our constraints
  and coherence are produced in a universal way, any reasonable data
  can be produced and any reasonable axioms will be satisfied.  For
  concreteness, we use the definition of~\cite{mccrudden:bal-coalgb}.

  The first piece of data we require to make $\cH(\lD)$ braided is a
  pseudo natural adjoint equivalence $\mathord{\otimes} \too[\eqv]
  \mathord{\otimes}\circ \tau$, where $\tau$ is the switch
  isomorphism.  This arises by \autoref{thm:h-locfr} from the braiding
  of \lD.  We also require two invertible modifications filling the
  usual hexagons for a braiding:
  \[\vcenter{\xymatrix@-1pc{
      & \mathclap{(A\ten B)\ten C}\phantom{C_C} \ar[dl]\ar[dr] \\
      (B\ten A)\ten C \ar[d] & \dtwocell\omit{\ze \iso} & A\ten (B\ten C)\ar[d]\\
      B\ten (A\ten C) \ar[dr] && (B\ten C)\ten A \ar[dl]\\
      & \mathclap{B\ten (C\ten A)}\phantom{C^C}
    }}\quad\text{and}\quad
  \vcenter{\xymatrix@-1pc{
      & \mathclap{A\ten (B\ten C)}\phantom{C_C} \ar[dl]\ar@{<-}[dr]\\
      A\ten (C\ten B)\ar[d] & \dtwocell\omit{\xi \iso}& (A\ten B)\ten C\ar[d]\\
      (A\ten C)\ten B \ar[dr] && C\ten (A\ten B) \ar@{<-}[dl]\\
      & \mathclap{(C\ten A)\ten B}\phantom{C^C}
    }}
  \]
  As before, since the corresponding hexagons commute in $\lD_0$, and
  by Lemmas \ref{thm:comp-compose} and \ref{thm:comp-ten} each side of
  each hexagon in $\cH(\lD)$ is a companion to the corresponding side
  in $\lD_0$, we have $\theta$ isomorphisms that we can take as $\ze$
  and $\xi$.  Finally, we must verify that the four 2-cell diagrams
  in~\cite[p136--139]{mccrudden:bal-coalgb} involving \ze\ and \xi\
  commute.  As with the axioms for a monoidal bicategory, both sides
  of these equalities are made up of $\theta$s relating companions of
  a single morphism in $\lD_0$, and thus by uniqueness they must be
  equal.



  Now suppose that \lD\ is symmetric.  To make $\cH(\lD)$ symmetric,
  we require first a \emph{syllepsis}, i.e.\ an invertible
  modification
  \[\vcenter{\xymatrix{A\ten B \ar@{=}[rr] \ar[dr] &
      \ar@{}[d]|-{\Downarrow \nu\iso}
      & A\ten B \ar@{<-}[dl]\\ &
      B\ten A }}\]
  Since the braiding in $\lD_0$ is self-inverse, the top and bottom of
  this triangle are both companions of $1_{A\ten B}$; thus we have a
  $\theta$ isomorphism between them which we take as $\nu$.  For
  $\cH(\lD)$ to be sylleptic, the syllepsis must satisfy the two
  axioms on~\cite[p144--145]{mccrudden:bal-coalgb}.  As before, these
  diagrams of 2-cells are made up entirely of $\theta$s relating
  companions of a single morphism in $\lD_0$, so they commute by
  uniqueness of $\theta$.



  Finally, for $\cH(\lD)$ to be symmetric, the syllepsis must satisfy
  one additional axiom, given on~\cite[p91]{mccrudden:bal-coalgb}.
  This follows automatically for the same reasons as before.
\end{proof}

Combining the arguments of Theorems \ref{thm:h-functor} and
\ref{thm:mon11-monbi}, we could show that passage from fibrant
monoidal double categories to monoidal bicategories is a functor of
tricategories, given a suitable definition of a tricategory of
monoidal bicategories.

\begin{rmk}
  Essentially the same proof as that of \autoref{thm:mon11-monbi}
  shows that any fibrant 2x1-category has an underlying tricategory.
  Note that unlike the construction of bicategories from
  1x1-categories (i.e.\ double categories), this case requires
  fibrancy even in the absence of monoidal structure, since the
  associativity and unit constraints of a 2x1-category are not 1-cells
  but rather morphisms of 0-cells.  There are many naturally occurring
  fibrant symmetric monoidal 2x1-categories, such as $\lD_0=$
  commutative rings, $\lD_1=$ algebras, and $\lD_2=$ modules, or the
  symmetric monoidal 2x1-category of \emph{conformal nets} defined
  in~\cite{bdh:confnets-i}.  All of these have underlying
  tricategories, which will be symmetric monoidal for any reasonable
  definition of symmetric monoidal tricategory.  More generally, as
  stated in \S\ref{sec:introduction}, we expect any fibrant $(n\times
  k)$-category to have an underlying $(n+k)$-category.
\end{rmk}

\bibliographystyle{alpha}
\bibliography{all,shulman}

\end{document}